\documentstyle[11pt]{article}
\newcommand{\ep}{\varepsilon}
\newcommand{\de}{\delta}
\newcommand{\G}{\Gamma}
\newcommand{\beq}[1]{\begin{equation}\label{#1}}
\newcommand{\eq}{\end{equation}}
\newcommand{\beqn}[1]{\begin{eqnarray}\label{#1}}
\newcommand{\eqn}{\end{eqnarray}}
\newcommand{\ioj}{J_0^{(j)}((1-q^2)zs;q^2)}
\newcommand{\iojm}{J_0^{(j)}((1-q^2)q^{-1}zs;q^2)}
\newcommand{\iojp}{J_0^{(j)}((1-q^2)qzs;q^2)}
\newcommand{\iijp}{J_1^{(j)}((1-q^2)q^{1-\frac\de2}zs;q^2)}
\newcommand{\iij}{J_1^{(j)}((1-q^2)q^{-\frac\de2}zs;q^2)}
\newcommand{\iojq}{ J_0^{(j)}((1-q^2)q^{1-\de}zs;q^2)}
\newcommand{\qpop}{\frac{1-q^2}2}
\newcommand{\pros}{\frac{2\partial_s}{1+q}}
\newcommand{\proz}{\frac{2\partial_z}{1+q}}
\newcommand{\p}{\partial}
\newcommand{\mbesj}{I_\nu^{(j)}((1-q^2)z;q^2)}

\newcommand{\mmbesj}{I_{-\nu}^{(j)}((1-q^2)z;q^2)}
\newcommand{\mbess}{I_\nu^{(2)}((1-q^2)z;q^2)}

\newcommand{\makdj}{K_\nu^{(j)}((1-q^2)z;q^2)}
\newcommand{\qqexp}{e_{q^2}(\frac{(1-q^2)^2}{4}z^2)}
\newcommand{\qqeexp}{E_{q^2}(-\frac{(1-q^2)^2}{4}q^2z^2)}
\newtheorem{predl}{Proposition}[section]
\newtheorem{defi}{Definition}[section]
\newtheorem{rem}{Remark}[section]
\newtheorem{lem}{Lemma}[section]
\newtheorem{cor}{Corollary}[section]

\begin{document}

\vspace{10mm}
\begin{flushright}
 ITEP-TH-06/01\\
\end{flushright}
\vspace{10mm}
\begin{center}
{\Large \bf The integral representations of the $q$-Bessel-Macdonald
functions}\\
\vspace{5mm}
V.-B.K.Rogov \footnote{The work was supported by the Russian Foundation for
Fundamental Research (grant no. 0001-00143) and the NIOKR MPS RF}\\
101475, Moscow, MIIT\\
e-mail vrogov@cemi.rssi.ru
\end{center}

\begin{abstract}

The $q$-Bessel-Macdonald functions of kinds 1, 2 and 3 are
considered. Their representations by classical integral are
constructed.
\end{abstract}

\section{Introduction}
\setcounter{equation}{0}

The definitions of the $q$-Bessel-Macdonald functions ($q$-BMF) and their
properties were given in \cite{OR1}. Their representations by the Jackson
$q$-integral have been constructed in \cite{OR2,R}.

The double integral appears necessarily if we consider the problems of the
harmonic analysis on the quantum Lobachevsky space. As the double Jackson
$q$-integral is connected hardly with the $q$-lattice we have not any
possibility to pass to the polar coordinates. So we are forced to introduce
the usual double integral for $q$-functions.

In the first sections we remind the known formulas for $q$-Bessel functions
and $q$-binomial.

\section{The modified $q$-Bessel functions and the $q$-BMF}
\setcounter{equation}{0}

\vspace{5mm}
In \cite{Ja} the $q$-Bessel functions were defined as follows:
\beq{2.1}
J_\nu^{(1)}(z,q)=\frac{(q^{\nu+1},q)_\infty}{(q,q)_\infty}(z/2)^\nu
\phantom1_2\Phi_1(0,0;q^{\nu+1};q,-\frac{z^2}{4}),
\eq
\beq{2.2}
J_\nu^{(2)}(z,q)=\frac{(q^{\nu+1},q)_\infty}{(q,q)_\infty}(z/2)^\nu
\phantom1_0\Phi_1(-;q^{\nu+1};q,-\frac{z^2q^{\nu+1}}{4}),
\eq
\beq{2.3}
J_\nu^{(3)}(z,q)=\frac{(q^{\nu+1},q)_\infty}{(q,q)_\infty}(z/2)^\nu
\phantom1_1\Phi_1(0;q^{\nu+1};q,-\frac{z^2q^{\frac{\nu+1}2}}{4}).
\eq
where $\phantom1_r\Phi_s$ is basic hypergeometric function \cite {GR},
$$
\phantom1_r\Phi_s(a_1,\cdots,a_r;b_1,\cdots,b_s;q,z)=
\sum_{n=0}^\infty\frac{(a_1,q)_n\ldots(a_r,q)_n}
{(q,q)_n(b_1,q)_n\ldots(b_s,q)_n}[(-1)^nq^{n(n-1)/2}]^{1+s-r}z^n,
$$
$$
(a,q)_n=\left\{
\begin{array}{lcl}
1&{\rm for}&n=0\\
(1-a)(1-aq)\ldots(1-aq^{n-1})&{\rm for}&n\ge1,\\
\end{array}
\right.
$$
$$
(a,q)_\infty=\lim_{n\to\infty}(a,q)_n,\qquad
(a_1,\ldots,a_k,q)_\infty=(a_1,q)_\infty\ldots(a_k,q)_\infty.
$$

It allows to introduce the modified $q$-Bessel functions
($q$-MBFs) using (\ref{2.1}), (\ref{2.2}) and (\ref{2.3}) similarly to
the classical case \cite{BE}.

\begin{defi}\label{d2.1}
The modified $q$-Bessel functions are the functions
$$
I_\nu^{(j)}(z,q)=\frac{(q^{\nu+1},q)_\infty}{(q,q)_\infty}(z/2)^\nu
\phantom1_\de\Phi_1\left(\underbrace{0,\ldots,0}_{\mbox{$\de$}};q^{\nu+1};
q,\frac{z^2q^{\frac{\nu+1}2(2-\de)}}{4}\right).
$$
Here
\beq{2.4}
\de=\left\{
\begin{array}{lcl}
2 & {\rm for} & j=1\\
0 & {\rm for} & j=2\\
1 & {\rm for} & j=3.\\
\end{array}
\right.
\eq
\end{defi}

Obviously,
$$
I_\nu^{(j)}(z,q)=e^{-\frac{i\nu\pi}{2}}J_\nu^{(j)}(e^{i\pi/2}z,q),
\qquad j=1,2,3.
$$
In the sequel we consider the functions
\beq{2.5}
I_\nu^{(1)}((1-q^2)z;q^2)=\sum_{k=0}^\infty\frac{(1-q^2)^k(z/2)^{\nu+2k}}
{(q^2,q^2)_k\Gamma_{q^2}(\nu+k+1)},\qquad |z|<\frac{2}{1-q^2},
\eq
\beq{2.6}
I_\nu^{(2)}((1-q^2)z;q^2)=
\sum_{k=0}^\infty\frac{q^{2k(\nu+k)}(1-q^2)^k(z/2)^{\nu+2k}}
{(q^2,q^2)_k\Gamma_{q^2}(\nu+k+1)}.
\eq
\beq{2.7}
I_\nu^{(3)}((1-q^2)z;q^2)=
\sum_{k=0}^\infty\frac{q^{k(\nu+k)}(1-q^2)^k(z/2)^{\nu+2k}}
{(q^2,q^2)_k\Gamma_{q^2}(\nu+k+1)},
\eq
where
$$
\Gamma_{q^2}(\nu)=\frac{(q^2,q^2)_{\infty}}{(q^{2\nu},q^2)_\infty}
(1-q^2)^{1-\nu}
$$
is the $q^2$-gamma function.
If $|q|<1,$ the series (\ref{2.6}) and (\ref{2.7}) are absolutely
convergent for all $z\ne0$. Consequently, $\mbess$ and
$I_\nu^{(3)}((1-q^2)z;q^2)$ are holomorphic functions outside a
neighborhood of zero.
\begin{rem}\label{r2.1}
$$
\lim_{q\to1-0}I_\nu^{(j)}((1-q^2)z;q^2)=I_\nu(z),\qquad j=1, 2, 3.
$$
\end{rem}
\begin{predl}\label{p2.1}
The function $\mbesj$ is a solution of the difference equation
\beq{2.8}
f(q^{-1}z)-(q^{-\nu}+q^\nu)f(z)+f(qz)=
q^{-\de}\frac{(1-q^2)^2}4z^2f(q^{1-\de}z),
\eq
where $j=1, 2, 3$ are connected with $\de=2, 0, 1$ by relations (\ref{2.4}).
\end{predl}
\begin{cor}\label{c2.1}
The function $I_{-\nu}^{(j)}((1-q^2)z;q^2)$ satisfies equation
(\ref{2.8}).
\end{cor}
\begin{defi}\label{d2.2}
We define the $q$-Bessel-Macdonald function ($q$-BMF) for $j=1, 2, 3$
as follows \cite{OR1,R}:
$$
\makdj=
$$
\beq{2.9}
=\frac12q^{-\nu^2+\nu}\G_{q^2}(\nu)\G_{q^2}(1-\nu)
\left[A_\nu^{|1-\de|}\mmbesj-A_{-\nu}^{|1-\de|}\mbesj\right],
\eq
where
\beq{2.10}
A_\nu=\sqrt{\frac{I_\nu^{(2)}(2;q^2)}{I_{-\nu}^{(2)}(2;q^2)}}.
\eq
\end{defi}
As in the classical case, this definition must be extended to integral
values of $\nu=n$ by passing to the limit in (\ref{2.9}).
\begin{defi}\label{d2.3}
The $q$-Wronskian of two solutions $f_\nu^1(z)$ and $f_\nu^2(z)$ of a
second-order difference equation is defined as follows:
$$
W(f_\nu^1,f_\nu^2)(z)=f_\nu^1(z)f_\nu^2(qz)-f_\nu^1(qz)f_\nu^2(z).
$$
\end{defi}

If the $q$-Wronskian does not vanish, then any solution of the
second-order difference equation can de written in form
$$
f_\nu(z)=C_1f_\nu^1(z)+C_2f_\nu^2(z).
$$
In this case the functions $f_\nu^1(z)$ and $f_\nu^2(z)$ form a fundamental
system of the solutions of the given equation.

\begin{predl}\label{p2.2}
The functions $\mbesj$ and $\makdj$ form a fundamental
system of the solutions of equation (\ref{2.8})
($z\ne\pm\frac{2q^{-r}}{1-q^2}, r=0, 1,\ldots$ if $j=1$).
\end{predl}

This Proposition is following from
$$
W(z)=\left\{
\begin{array}{lcl}
\frac{q^{-\nu}(1-q^2)}2A_\nu\qqexp & {\rm for}&\de=2\\
\frac{q^{-\nu}(1-q^2)}2 & {\rm for}&\de=1\\
\frac{q^{-\nu}(1-q^2)}2A_\nu\qqeexp &{\rm for}&\de=0.\\
\end{array}
\right.
$$
(See \cite{R}.)
Obviously, this function is defined for $z\ne\pm\frac{2q^{-r}}{1-q^2},
r=0, 1,\ldots $ if $j=1 ~~(\de=2)$ and does not vanish.

\section{Some ancillary formulas}
\setcounter{equation}{0}

There is a $q$-analog of the classical binomial formula \cite{GR}
$$
(1-z)^{-a}=\sum_{k=0}^\infty\frac{(a)_k}{k!}z^k,~(a)_k=
\frac{\Gamma(a+k)}{\Gamma(a)},
\qquad |z|<1,
$$
$$
\frac{(q^\alpha z,q)_\infty}{(z,q)_\infty}=\sum_{k=0}^\infty
\frac{(q^\alpha,q)_k}{(q,q)_k}z^k,\qquad |z|<1.
$$
We need in two generalizations of the $q$-binomial
\begin{equation}\label{3.1}
r(a,b,z,q)=\frac{(az,q)_\infty}{(bz,q)_\infty}
\end{equation}
\begin{equation}\label{3.2}
R(a,b,\gamma,z,q^2)=\frac{(az^2,q^2)_\infty}
{(bz^2,q^2)_\infty}z^\gamma
\end{equation}
\begin{predl}\label{p3.1}
The function $R(a,b,\gamma,z,q^2)$ (\ref{3.2}) satisfies
the difference equation
\begin{equation}\label{3.3}
z^2[bq^\gamma R(a,b,\gamma,z,q^2)-aR(a,b,\gamma,qz,q^2)]=
q^\gamma R(a,b,\gamma,z,q^2)-R(a,b,\gamma,qz,q^2).
\end{equation}
\end{predl}
The {\bf Proof} see in \cite{OR2}

It was shown in \cite{OR2} that if $\alpha>\beta$ then
\beq{3.4}
\frac{(-q^{2\alpha}z^2,q^2)_\infty}{(-q^{2\beta}z^2,q^2)_\infty}=
\frac{(q^{2(\alpha-\beta)},q^2)_\infty}{(q^2,q^2)_\infty}
\sum_{k=0}^\infty\frac{(q^{2(\beta-\alpha+1)},q^2)_kq^{2(\alpha-\beta-1)k}}
{(q^2,q^2)_k(1+z^2q^{2\beta+2k})}.
\eq

\begin{rem}\label{r3.1}
Let $a=\epsilon q^{2\alpha}, b=\epsilon q^{2\beta}, \epsilon=\pm1,$
in (\ref{3.2}). Then if $q\to1-0$ the difference equation (\ref{3.3})
takes the form of the differential equation
\beq{3.5}
z(1-\epsilon z^2)R'(z)-[\gamma+\epsilon(2\alpha-2\beta-\gamma)z^2]R(z)
=0
\eq
with solution
$$
R(z)=Cz^\gamma(1-\epsilon z^2)^{\beta-\alpha}.
$$
\end{rem}

Designate the $q$-derivative of $g(x)$ by
$$
\p_xg(x)=\frac{g(x)-g(qx)}{(1-q)x}.
$$
It follows from
\beq{3.6}
\ioj=\sum_{k=0}^\infty\frac{(-1)^kq^{(2-\de)k^2}(1-q^2)^{2k}(sz/2)^{2k}}
{(q^2,q^2)_k^2}
\eq
that
\beq{3.7}
\proz\iojm=-q^{1-\frac\de2}s\iij,
\eq
\beq{3.8}
\proz\ioj=-q^{1-\frac\de2}s\iijp
\eq
and
\beq{3.9}
\pros s\iijp=q^{-\frac\de2}zs\iojq
\eq

\begin{lem}\label{l3.1}
If $F(x)$ is a differentiable (in classical sense) function in some
neighborhood of zero then
\beq{3.10}
\lim_{\ep\to0}\int_{q\ep}^\ep\frac{F(x)}xdx=-F(0)\ln q.
\eq
\end{lem}
{\bf Proof.} Integrating by parts we have
$$
\lim_{\ep\to0}\int_{q\ep}^\ep\frac{F(x)}xdx=
\lim_{\ep\to0}\left[F(x)\ln x|_{q\ep}^\ep-
\int_{q\ep}^\ep F'(x)\ln xdx\right]=
$$
$$
=\lim_{\ep\to0}\left[(F(\ep)-F(q\ep))\ln\ep-F(q\ep)\ln q-
\int_{q\ep}^\ep F'(x)\ln xdx\right].
$$
Using Lagrange's theorem and the theorem about mean value we obtain in the
left side
$$
\lim_{\ep\to0}\left[F'(\theta_1\ep)(1-q)\ep\ln\ep-F(q\ep)\ln q-
(1-q)\ep F'(\theta_2\ep)\ln(\theta_2\ep)\right],
$$
where $\theta_1\in (0,1), ~\theta_2\in (0,1)$.
As $\lim_{\ep\to0}\ep\ln\ep=0$ we have (\ref{3.10}). \rule{5pt}{5pt}

\begin{lem}\label{l3.2}
If $f(x)$ and $g(x)$ are integrable on $(0,\infty)$ and differentiable ones in
zero (in classical sense) then the $q$-analog of formula of integration by
parts for the $q$-derivative takes place
\beq{3.11}
\int_0^\infty\p_x f(x)g(x)dx=f(0)g(0)\frac{\ln q}{1-q}-
\int_0^\infty f(qx)\p_xg(x)dx.
\eq
\end{lem}
{\bf Proof.}
$$
\int_0^\infty\p_x f(x)g(x)dx=\int_0^\infty\frac{f(x)-f(qx)}{(1-q)x}g(x)dx=
$$
$$
=\int_0^\infty\frac{f(x)g(x)-f(qx)g(qx)}{(1-q)x}dx-
\int_0^\infty f(qx)\frac{g(x)-g(qx)}{(1-q)x}dx=
$$
$$
=\int_0^\infty\p_x(f(x)g(x))dx-\int_0^\infty f(qx)\p_xg(x)dx.
$$
Using Lemma \ref{l3.1} calculate the first integral in the right side.
$$
\int_0^\infty\p_x(f(x)g(x))dx=
\frac1{1-q}\lim_{\ep\to0}\left[\int_\ep^\infty\frac{f(x)g(x)}xdx-
\int_\ep^\infty\frac{f(qx)g(qx)}xdx\right]=
$$
$$
=-\frac1{1-q}\lim_{\ep\to0}\int_{q\ep}^\ep\frac{f(x)g(x)}xdx=
\frac{\ln q}{1-q}f(0)g(0).
$$
The statement of the Lemma follows from here. \rule{5pt}{5pt}

\begin{cor}\label{c3.1}
If  $f(s)g(s)s^{-1}$ is integrable function on $(0,\infty)$ then
\beq{3.12}
\int_0^\infty f(s)\p_sg(s)ds=-\int_0^\infty\p_sf(s)g(qs)ds.
\eq
\end{cor}
The {\bf Proof} follows from $f(0)g(0)=0$ in this case. \rule{5pt}{5pt}

\section{The integral representations of the $q$-BMFs}
\setcounter{equation}{0}

\vspace{5mm}
We will assume that $z$ and $s$ are the commuting variables.
\begin{predl}\label{p4.1}
$q$-BMF $\makdj$ for ${\rm Re}\nu>0$ can be represented as the integral
$$
\makdj=-\frac{q^{-\nu^2+\nu(1-\de)}(1-q^2)}{2\ln q}\G_{q^2}(\nu+1)
A_\nu^{|1-\de|}\times
$$
\beq{4.1}
\times(z/2)^{-\nu}\int_0^\infty\frac{(-q^{2\nu+2-\de\nu}s^2,q^2)_\infty}
{(-q^{-\de\nu}s^2,q^2)_\infty}s\ioj ds,
\eq
where the constant $A_\nu$ is defined by (\ref{2.10}) and $j=1, 2, 3$ are
connected with $\de=2, 0, 1$ by relations (\ref{2.4}).
\end{predl}

{\bf Proof.} Consider the absolutely convergent integral
\beq{4.2}
S^{(j)}(z)=\int_0^\infty f_\nu^{(j)}(s)sJ_0^{(j)}((1-q^2)zs;q^2)ds,
\eq
and require that $S^{(j)}(z)(z/2)^{-\nu}$ satisfies the difference equation
(\ref{2.8}). Then $S^{(j)}(z)$ satisfies the equation
\beq{4.3}
S^{(j)}(q^{-1}z)-S^{(j)}(z)-q^{-2\nu}[S^{(j)}(z)-S^{(j)}(qz)]=
q^{\nu(\de-2)-\de}(\qpop)^2z^2S^{(j)}(q^{1-\de}z).
\eq
Substituting (\ref{4.2}) in (\ref{4.3}) and multiplying it on
$\frac{2z^{-1}}{1-q^2}$ we obtain
$$
\int_0^\infty f_\nu^{(j)}(s)s\frac{2z^{-1}}{1-q^2}[\iojm-\ioj]ds-
$$
$$
-q^{2\nu}\int_0^\infty f_\nu^{(j)}(s)s\frac{2z^{-1}}{1-q^2}[\ioj-\iojp]ds=
$$
$$
=q^{\nu(\de-2)-\de}\qpop z\int_0^\infty f_\nu^{(j)}(s)s\iojq ds.
$$

Due to (\ref{3.7}) - (\ref{3.9}) we can write
$$
\int_0^\infty f_\nu^{(j)}(s)s^2\iij ds-
q^{-2\nu+1}\int_0^\infty f_\nu^{(j)}(s)s^2\iijp ds=
$$
$$
=q^{\nu(\de-2)-\frac\de2}\qpop s\int_0^\infty f_\nu^{(j)}(s)s\iojq ds
$$
or
$$
\int_0^\infty f_\nu^{(j)}(s)s^{2\nu+2}\pros[s^{-2\nu+1}\iij]ds=
$$
$$
=-q^{\nu(\de-2)}\int_0^\infty f_\nu^{(j)}(s)\pros\iij ds.
$$

Using (\ref{3.12}) we obtain
$$
\int_0^\infty\p_s(f_\nu^{(j)}(s)s^{2\nu+2})q^{-2\nu+1}s^{-2\nu+1}\iijp ds=
$$
$$
=-q^{\nu(\de-2)}\int_0^\infty\p_sf_\nu^{(j)}(s)qs\iijp ds.
$$
Thus we come to the difference equation for $f_\nu^{(j)}(s)$
\beq{4.4}
q^{-\de\nu}s^2[-f_\nu^{(j)}(s)+q^{2\nu+2}f_\nu^{(j)}(qs)]=
f_\nu^{(j)}(s)-f_\nu^{(j)}(qs).
\eq
It follows from Proposition \ref{p3.1} the function
\beq{4.5}
f_\nu^{(j)}(s)=\frac{(-q^{2\nu+2-\de\nu}s^2,q^2)_\infty}
{(-q^{-\de\nu}s^2,q^2)_\infty}
\eq
satisfies to (\ref{4.4}), and it follows from (\ref{3.4}) integral
(\ref{4.2}) is absolutely convergent.

As $S^{(j)}(z)(z/2)^{-\nu}$ is a solution to  (\ref{2.8}) it can be
represented as
\beq{4.6}
S^{(j)}(z)(z/2)^{-\nu}=A\mbesj+B\makdj.
\eq

Let $j=1$. As it follows from \cite{OR1} $I_\nu^{(1)}((1-q^2)z;q^2)$ is a
meromorphic function with the ordinary poles $z=\pm\frac{2q^{-r}}{1-q^2},
~r=0,1,\ldots$, and $K_\nu^{(1)}((1-q^2)z;q^2)$ and the left side of
(\ref{4.6}) are the holomorphic functions in region ${\rm Re}z>0$.

Let $j=2, 3$. It is easily to show (see \cite{R}) that $\lim_{z\to\infty}
\mbesj=\infty,\\ \lim_{z\to\infty}\makdj=0$ and the left side of (\ref{4.6})
tends to zero if $z\to\infty$.

So for any $j=1, 2, 3 ~~A=0$. Multiplying the both sides of (\ref{4.6}) on
$(z/2)^\nu$ and putting $z=0$ we obtain from (\ref{2.9}) and (\ref{2.5}) - (\ref{2.7})
\beq{4.7}
\int_0^\infty\frac{(-q^{2\nu+2-\de\nu}s^2,q^2)_\infty}
{(-q^{-\de\nu}s^2,q^2)_\infty}sds=\frac{B}2q^{-\nu^2+\nu}\G_{q^2}(\nu)
A_\nu^{|1-\de|}.
\eq

Now calculate the integral in the left side of (\ref{4.7}).
$$
\int_0^\infty\frac{(-q^{2\nu+2-\de\nu}s^2,q^2)_\infty}
{(-q^{-\de\nu}s^2,q^2)_\infty}sds=\frac{q^{\de\nu}}2
\int_0^\infty\frac{(-q_1^{\nu+1}x,q_1)_\infty}
{(-x,q_1)_\infty}dx=
\frac{q^{\de\nu}}2\int_0^\infty e_{q_1}(-x)E_{q_1}(q_1^{\nu+1}x)dx,
$$
where $x=q^{-\de\nu}s^2$ and $q_1=q^2$ \cite{GR}.

Note, that
$$
\p_xe_{q_1}(-x)=-\frac1{1-q_1}e_{q_1}(-x),
~~~~\p_xE_{q_1}(q_1^{\nu+1}x)=\frac{q_1^{\nu+1}}{1-q_1}E_{q_1}(q_1^{\nu+2}x).
$$
Using Lemma \ref{l3.2} we obtain
$$
\int_0^\infty e_{q_1}(-x)E_{q_1}(q_1^{\nu+1}x)dx=
-\ln q_1 +q_1^\nu\int_0^\infty e_{q_1}(-x)E_{q_1}(q_1^{\nu+1}x)dx.
$$
So
$$
\int_0^\infty e_{q_1}(-x)E_{q_1}(q_1^{\nu+1}x)dx=-\frac{\ln q_1}{1-q_1^\nu},
$$
and we have
\beq{4.8}
\int_0^\infty\frac{(-q^{2\nu+2-\de\nu}s^2,q^2)_\infty}
{(-q^{-\de\nu}s^2,q^2)_\infty}sds=-\frac{q^{\de\nu}\ln q}{1-q^{2\nu}}.
\eq
It follows from (\ref{4.8}) that
$$
B=-\frac{2q^{\nu^2-\nu(1-\de)}\ln q}{(1-q^2)\G_{q^2}(\nu+1)A_\nu^{|1-\de|}},
$$
and we have (\ref{4.1}).\rule{5pt}{5pt}

\begin{rem}\label{r4.1}
It follows from \cite[Remark 5.1]{OR1} $A_\nu \to 1$ if $q \to 1-0$,
and it follows from Remark \ref{r3.1} that if $q \to 1-0$ we come to the
classical integral representation of Bessel-Macdonald function \cite{BE}
$$
K_\nu(z)=\G(\nu+1)\left(\frac z2\right)^{-\nu}
\int_0^\infty(1+s^2)^{-\nu-1}sJ_0(zs)ds.
$$
\end{rem}

\section{The representation of the $q$-BMFs by a double integral}
\setcounter{equation}{0}

We take $z, ~s\in\mathbf{C}$ in this section.

Consider function
\beq{5.1}
\xi_\eta^{(\de)}(s)=\sum_{n=0}^\infty\frac{q^{(2-\de)\eta n^2}(1-q^2)^n}
{(q^2,q^2)_n}s^n, ~~~\eta\ge0.
\eq
Obviously this series converges for any $s$ if $(2-\de)\eta>0$ and for
$|s|<\frac1{1-q^2}$ if $(2-\de)\eta=0$. Consider three cases.\\
1. $\eta=0.$ It follows from (\ref{5.1})
\beq{5.2}
\xi_0^{(\de)}(s)=e_{q^2}((1-q^2)s), ~~~\de=2,0,1.
\eq
2. $\eta=1.$
\beq{5.3}
\xi_1^{(\de)}(s)=\left\{
\begin{array}{lcl}
e_{q^2}((1-q^2)s) & {\rm for} & \de=2\\
E_{q^2}((1-q^2)qs) & {\rm for} & \de=0\\
\phantom._1\Phi_1(0;-q;q,-(1-q^2)q^{\frac12}s) & {\rm for} & \de=1.\\
\end{array}
\right.
\eq
3. $\eta=\frac12$. It follows from (\ref{5.1})
\beq{5.4}
\xi_{\frac12}^{(\de)}(s)=\left\{
\begin{array}{lcl}
e_{q^2}((1-q^2)s) & {\rm for} & \de=2\\
\phantom._1\Phi_1(0;-q;q,-(1-q^2)q^{\frac12}s) & {\rm for} & \de=0.\\
\phantom._3\Phi_3(0,0,0;-q^{\frac12},iq^{\frac12},-iq^{\frac12};
q^{\frac12},-(1-q^2)q^{\frac14}s) & {\rm for} & \de=1.\\
\end{array}
\right.
\eq

Assume $s=\rho e^{i\phi}, ~z=re^{i\psi}$, and consider integral
\beq{5.5}
J(r\rho)=\int_{-\pi}^{\pi}\xi_\eta^{(\de)}(ir\rho e^{-i(\psi+\phi)})
\xi_{1-\eta}^{(\de)}(ir\rho e^{i(\psi+\phi)})d\phi.
\eq
Let $\de<2, ~\eta=\frac12$. In this case we can calculate this integral
term by term.
$$
J(r\rho)=
\int_{-\pi}^{\pi}\sum_{n=0}^\infty\frac{q^{(2-\de)\frac{n^2}2}(1-q^2)^n}
{(q^2,q^2)_n}(ir\rho)^ne^{-in\phi}
\sum_{m=0}^\infty\frac{q^{(2-\de)\frac{m^2}2}(1-q^2)^m}
{(q^2,q^2)_m}(ir\rho)^me^{im\phi}d\phi=
$$
$$
=\sum_{n=0}^\infty\frac{q^{(2-\de)\frac{n^2}2}(1-q^2)^n}
{(q^2,q^2)_n}(ir\rho)^n
\sum_{m=0}^\infty\frac{q^{(2-\de)\frac{m^2}2}(1-q^2)^m}
{(q^2,q^2)_m}(ir\rho)^m\int_0^{2\pi}e^{i(m-n)\phi}d\phi=
$$
$$
=2\pi\sum_{n=0}^\infty(-1)^n\frac{q^{(2-\de)n^2}(1-q^2)^{2n}}
{(q^2,q^2)_n^2}(r\rho)^{2n}.
$$
It follows from (\ref{3.6}) the last series is $q^2$-Bessel function and
we have
\beq{5.6}
J_0^{(j)}((1-q^2)2r\rho;q^2)=
\frac1{2\pi}\int_{-\pi}^\pi\xi_{\frac12}^{(\de)}(ir\rho e^{-i(\psi+\phi)})
\xi_{\frac12}^{(\de)}(ir\rho e^{i(\psi+\phi)})d\phi,
\eq
where $\xi_{\frac12}^{(\de)}$ is defined by (\ref{5.4}) and $j=2,3$ are connected
with $\de=0,1$ by relations (\ref{2.4}).

Let $\de<2$ and $\eta=0$. Then we have the same result for
$r\rho<\frac1{1-q^2}$, e.i.
\beq{5.7}
J_0^{(j)}((1-q^2)2r\rho;q^2)=\frac1{2\pi}\int_{-\pi}^\pi
e_{q^2}(i(1-q^2)r\rho e^{-i(\psi+\phi)})
\xi_1^{(\de})(ir\rho e^{i(\psi+\phi)}) d\phi,
\eq
where $\xi_1^{(de)}$ is defined by (\ref{5.3}).
The left side is a holomorphic function outside a neighborhood of zero,
and so we can consider $J_0^{(j)}((1-q^2)2r\rho;q^2), ~j=2,3$ as the
analytic continuation of (\ref{5.8}).

Let $\de=2$.  It follows from (\ref{5.2}) - (\ref{5.6}) for
$r\rho<\frac1{1-q^2}$
\beq{5.8}
J_0^{(1)}((1-q^2)2r\rho;q^2)=\frac1{2\pi}\int_{-\pi}^\pi
e_{q^2}(i(1-q^2)r\rho e^{-i(\psi+\phi)})
e_{q^2}(i(1-q^2)r\rho e^{i(\psi+\phi)})d\phi.
\eq
In this case the both sides of (\ref{5.8}) are the meromorphic functions
with the ordinary poles in points $r\rho=\pm\frac{iq^{-k}}{1-q^2}$. So we will
consider $J_0^{(1)}((1-q^2)2r\rho;q^2)$ as a analytic continuation of
(\ref{5.8}).

Now we can formulate
\begin{predl}\label{p5.1}
The $q$-BMF can be represented by double integral
\beq{5.9}
K\nu^{(j)}(2(1-q^2)|z|,q^2)=
-\frac{q^{-\nu^2+\nu(1-\de)}(1-q^2)}{8\pi\ln q}\G_{q^2}(\nu+1)
A_\nu^{|1-\de|}|z|^{-\nu}\times
\eq
$$
\times\int\int\frac{(-q^{2\nu+2-\de\nu}\bar ss,q^2)_\infty}
{(-q^{-\de\nu}\bar ss,q^2)_\infty}e_{q^2}(i(1-q^2)\bar z\bar s)
\xi_1^{(\de)}(izs)d\bar sds
$$
or
\beq{5.10}
K\nu^{(j)}(2(1-q^2)|z|,q^2)=
-\frac{q^{-\nu^2+\nu(1-\de)}(1-q^2)}{8\pi\ln q}\G_{q^2}(\nu+1)
A_\nu^{|1-\de|}|z|^{-\nu}\times
\eq
$$
\times\int\int\frac{(-q^{2\nu+2-\de\nu}\bar ss,q^2)_\infty}
{(-q^{-\de\nu}\bar ss,q^2)_\infty}\xi_{\frac12}^{(\de)}(i\bar z\bar s)
\xi_{\frac12}^{(\de)}(izs)d\bar sds,
$$
where $\xi_\eta^{(\de)}$ are defined by (\ref{5.2}) - (\ref{5.4}),
the constant $A_\nu$ is defined by (\ref{2.10}) and $j=1, 2, 3$ are
connected with $\de=2, 0, 1$ by relations (\ref{2.4}).
\end{predl}
{\bf Proof.} Substituting (\ref{5.6}) - (\ref{5.8}) in (\ref{4.1}) we
obtain (\ref{5.9}) or (\ref{5.10}) respectively. \rule{5pt}{5pt}

\begin{rem}\label{r5.1}
It ie easily to show that if $q \to 1-0$ the all functions (\ref{5.2}) -
(\ref{5.4}) tend to the usual exponential $e^s$, and we come to the
classical integral representation of Bessel-Macdonald function
$$
K_\nu(2\sqrt{\bar zz})=\G(\nu+1)(\sqrt{\bar zz})^{-\nu}
\int\int(1+\bar ss)^{-\nu-1}\exp(i(\bar z\bar s+zs))d\bar sds.
$$
\end{rem}

\small{
 }

\end{document}